\documentclass[11pt]{article}
\usepackage{latexsym}
\usepackage{amsfonts,amssymb,amsmath}

\newtheorem{thm}{Theorem}[section]

\newtheorem{lem}[thm]{Lemma}

\newtheorem{defn}[thm]{Definition}

\newcommand{\ov }{\overline }

\title{New embeddings between the Higman-Thompson groups }

\author{ J.C.\ Birget }

\date{\scriptsize{
10.iii.2019}}

\begin{document}
\maketitle

\begin{abstract}
We give a direct proof that all Higman-Thompson groups of the form $G_{k,1}$
(for $k \ge 2$) are embedded in one another, which is a recent result of 
N.\ Matte Bon. This extends the embeddings given by Higman in 1974.
\end{abstract}

\section{Introduction} 

Higman \cite[Theorem 7.2 and Lemma 7.1]{Hig74} showed that if
$\,K = 1 + (k-1) \, d \,$ for some $d \ge 1$, then $G_{K,r} \le G_{k,r}$.
In particular, $G_{K,1} \le G_{2,1}$ for all $K \ge 2$. 
It is also known that all $G_{k,1}$ are non-isomorphic for different 
$k\,$ (\cite[Theorem 6.4]{Hig74}, later generalized to the groups $G_{k,r}$
by \cite{Pardo}).

For a long time there has been a common belief that it was well known, and 
easy to prove, 
that all $G_{k,1}$ embed into each other (see the comments in section 3). 
However, the first proof of this is quite recent; it follows from 
\cite[Coroll.\ 11.16]{MatteBon} near the end of a long paper by Nicol\'as 
Matte Bon.

\begin{thm} \label{THMembed2intok} {\rm (N.\ Matte Bon).}
 \ All Higman-Thompson groups $G_{k,1}$ (for $k \ge 2$) are embedded in one 
another; i.e., for all $i, j \ge 2:$ \ $G_{i,1} \le G_{j,1}$. 
\end{thm}

The contribution of the present paper is a proof that is direct, elementary, 
and relatively short.
In section 2 we prove that $G_{2,1} \le G_{k,1}$ for all $k > 2$; this, in 
combination with Higman's embeddings, implies the Theorem.

\bigskip

The Higman-Thompson groups $G_{k,r}$ (for $k \ge 2$, $k > r \ge 1$) were 
introduced by Graham Higman \cite{Hig74} as a generalization of the 
Thompson group $V$ ($= G_{2,1}$) \ \cite{Th0, Th}.  
We refer to the literature (in particular \cite{Th0, McKTh, Th, Hig74, 
Scott, CFP}) for some of the remarkable properties of these groups; these
groups occur in many subjects (e.g., Pardo used connections with Leavitt 
path algebras to prove his result \cite{Pardo}). 

To define $G_{k,1}$ we follow \cite{BiThomps} (which is similar to 
\cite{Scott} except for terminology).
We use the alphabet $A_k = \{a_0, a_1, \, \ldots \, , a_{k-1}\}$, for any 
integer $k \ge 2$. Often we just write $A$ instead of $A_k$.
The empty string is denoted by $\varepsilon$, and the set of all strings 
over $A$ is denoted by $A^*$, and the set of all non-empty strings is 
denoted by $A^+$.  For a string $x \in A^*$, $|x|$ denotes the length.  
For a set $S \subseteq A^*$, $|S|$ denotes the cardinality.  
Concatenation of sets $S, T \subseteq A^*$ is denoted by $ST$ or 
$S \cdot T$, and defined by $ST = \{st : s \in S, \, t \in T\}$.
For $x, p \in A^*$ we say that $p$ is a {\it prefix} of $x$ iff 
$(\exists u \in A^*) \, x = pu$; this is denoted by $p \le_{\rm pref} x$. 
Two strings $x, y \in A^*$ are called {\it prefix-comparable} (denoted 
by $x \, \|_{\rm pref} \, y \,$) iff 
$\, x \le_{\rm pref} y \,$ or $ \, y \le_{\rm pref} x$. 
A {\it prefix code} is any subset $P \subset A^*$ such that for all 
$p_1, p_2 \in P$: $\, p_1 \, \|_{\rm pref} \, p_2 \,$ implies $p_1 = p_2$.
A {\it right ideal} of $A^*$ is any subset $R \subseteq A^*$ such that 
$R = R \cdot A^*$. A subset $C \subseteq R$ generates $R$ as a right ideal 
iff $R = C \cdot A^*$. 
It is easy to prove that every finitely generated right ideal is generated 
by a unique finite prefix code, and this prefix code is the minimum 
generating set of the right ideal (with respect to $\subseteq$). 
A {\it maximal prefix code} is a prefix code $P \subset A^*$ that is not a 
strict subset of any other prefix code of $A^*$. 

In this paper, {\em function} means partial function. For a function 
$f: A^* \to A^*$, the domain and image sets are denoted by ${\rm Dom}(f)$,
respectively ${\rm Im}(f)$.
A {\it right ideal morphism} of $A^*$ is a function $f: A^* \to A^*$ such 
that for all $x \in {\rm Dom}(f)$ and all $w \in A^*$:

\smallskip

 \ \ \ $\, f(xw) = f(x) \ w$.

\smallskip

\noindent
In that case, ${\rm Dom}(f)$ is a right ideal; one easily proves that 
${\rm Im}(f)$ is also a right ideal. The prefix code that generates
${\rm Dom}(f)$ is denoted by ${\rm domC}(f)$, and is called the {\it domain
code} of $f$; the prefix code that generates ${\rm Im}(f)$ is denoted by 
${\rm imC}(f)$, and is called the {\it image code}.
The following inverse monoid is a stepping stone towards defining $G_{k,1}$:

\medskip
 
\begin{minipage}{\textwidth}
${\cal RI}_A^{\sf fin}$ $\, = \,$  $\{ f : f$ is a right ideal morphism of 
    $A^*$, $\, f$ is {\em injective}, and  

\hspace{0.9in}
    ${\rm domC}(f)$ and ${\rm imC}(f)$ are {\em finite maximal} 
    prefix codes\}.
\end{minipage}

\medskip

\noindent We also write ${\cal RI}^{\sf fin}_k$ for ${\cal RI}_A^{\sf fin}$ 
(where $k = |A|$).  It is proved in \cite[Prop.\ 2.1]{BiThomps} that every 
$f \in {\cal RI}^{\sf fin}_k$ is contained in a unique $\subseteq$-maximum 
right ideal morphism in ${\cal RI}^{\sf fin}_k$; this is called the 
{\it maximum extension} of $f$.  The Higman-Thompson group $G_{k,1}$ (where 
$k = |A|$) is a homomorphic image of ${\cal RI}^{\sf fin}_k$, and also a 
subset of ${\cal RI}^{\sf fin}_k$ (as a set):

\begin{defn} \label{ThompsV} {\bf (the Higman-Thompson group $G_{k,1}$).}
 \ The Higman-Thompson group $G_{k,1}$, as a set, consists of the right 
ideal morphisms $f \in {\cal RI}_k^{\sf fin}$ that are maximum extensions 
in ${\cal RI}_k^{\sf fin}$. 
The multiplication in $G_{k,1}$ consists of composition, followed by maximum
extension.
\end{defn}
See \cite{BiThomps} for a proof that this multiplication turns $G_{k,1}$
into a group.

Every element $f \in {\cal RI}_k^{\sf fin}$ (and in particular, every 
$f \in G_{k,1}$) is determined by the restriction of $f$ to ${\rm domC}(f)$; 
this is a bijection from the finite prefix code ${\rm domC}(f)$ onto the 
finite prefix code ${\rm imC}(f)$. We call such a finite bijection a 
{\it table}. 
We do not assume here that $f$ is a maximum extension, so for an element 
$f \in G_{k,1}$ there are many non-maximal tables that determine $f$ by
maximal extension.
The well-known tree representation of $G_{k,1}$ is obtained by using the 
prefix trees of ${\rm domC}(f)$ and ${\rm imC}(f)$.

\begin{lem} \label{extStep} 
 \ The right ideal morphism $f \in {\cal RI}_A^{\sf fin}$ determined by a 
table $F$: $P \to Q$ can be extended \ iff \ there exist $p, q \in A^*$ such 
that for every $\alpha \in A$:
 \ $p \alpha \in P$, \ $q \alpha \in Q$, and $\, F(p \alpha) = q \alpha$.

In that case, $f$ can be extended by defining $\, f(p) = q$. So the table 
for this extension is obtained be replacing 
$\, \{(p\alpha, q\alpha) : \alpha \in A\}\,$ by $\{(p,q)\}$ in the table. 
This is called an {\em extension step} of the table $F$.  
\end{lem}
{\sc Proof.} See \cite[Lemma 2.2]{BiThomps}. \ \ \  \ \ \  $\Box$

\medskip

\noindent
Since in an extension step the cardinality of ${\rm domC}(f)$ decreases,
only finitely steps are needed to reach the maximum extension of $f$.

\bigskip

\noindent {\bf Notation.} For any prefix code $P \subseteq \{a_0,a_1\}^*$,
 \ ${\sf spref}(P)$  denotes the set of {\em strict prefixes} of the 
elements of $P$.  Formally, 
 
\smallskip

 \ \ \  \ \ \ ${\sf spref}(P) \, = \, $
$\{x \in \{a_0,a_1\}^* : \ (\exists p \in P)\,[\,x <_{\rm pref} p \,]\,\}$.

\smallskip

\begin{lem} \label{MaxprefCodeAk} \hspace{-0.07in}.

\noindent {\small \rm (1)} \ If $P$ is a finite maximal prefix code over 
$A_2 = \{a_0, a_1\}$,  then 
$ \ P \ \cup \ {\sf spref}(P) \cdot \{a_2, \,\ldots\, , a_{k-1}\} \ $  
is a finite maximal prefix code over 
$A_k = \{a_0, a_1, \, \ldots \, , a_{k-1}\}$.

\smallskip

\noindent {\small \rm (2)} \ If 
$ \ P \ \cup \ Q \cdot \{a_2, \, \ldots \, , a_{k-1}\} \,$
is a finite maximal prefix code over $A_k$, where $P$ and  $Q$ are finite 
subsets of $\{a_0,a_1\}^*$, then $P$ is a finite maximal prefix code over 
$\{a_0, a_1\}$ and $\, Q = {\sf spref}(P)$.
\end{lem}
{\sc Proof.} See \cite[Lemma 9.1]{BiCoNP}, where $k = 3$; the general case 
is similar.  
 \ \ \ $\Box$

\begin{lem} \label{OmegaPref}
 \ Let $P \subset A^*$ be a finite maximal prefix code. Then every
$v \in A^{\omega}$ has a unique prefix in $P$. \ Formally,
 \ $(\forall v \in A^{\omega})(\exists! \,p \in P,$
$u \in A^{\omega}) \, [\, v = p u \,]$.
\end{lem}
{\sc Proof.} The proof is straightforward. 
 \ \ \ $\Box$

\section{Embedding {\boldmath $G_{2,1}$} into {\boldmath $G_{k,1}$} }

\noindent To embed $G_{2,1}$ into $G_{k,1}$ the following subgroup of 
$G_{k,1}$ is used as an intermediary stage:

\bigskip

\begin{minipage}{\textwidth}
 \ \ \ $G_{k,1}(0,1|2|\ldots|k{\rm -}1)$ 

\smallskip

\hspace{0.22in} $ = \ \big\{g \in G_{k,1} : \ \ $  
   {\small (1)} \ \ ${\rm domC}(g) \,\cup\, {\rm imC}(g)$  
   $ \ \subset \ $  $\{a_0,a_1\}^* \ \cup \ $
   $\bigcup_{i=1}^{k-1} \{a_0,a_1\}^* a_i$,
  \ and

\hspace{1.3in} {\small (2)} \ \ for all $x \in A_k^{\,*}$:

\hspace{1.3in} {\small (2.1)}
 \ \ \ \ \ $x \in \{a_0,a_1\}^*$ $ \ \Leftrightarrow \ $ 
 $g(x) \in \{a_0,a_1\}^*$, \ and

\hspace{1.3in} {\small (2.2)}
 \ \ \ \ \ for all $\,i = 2, \ldots, k{\rm -}1: \ \ x \in \{a_0,a_1\}^*a_i$
$ \ \Leftrightarrow \ $
$g(x) \in \{a_0,a_1\}^* a_i \ \big\}$.
\end{minipage}

\bigskip

\noindent The special case $G_{3,1}(0,1|2)$ was introduced in 
\cite[Def.\ 4.4]{BiCoNP} (in \cite{BiCoNP} the alphabet 
$A_3 = \{a_0,a_1,a_2\}$ was denoted by $\{0,1,\#\}$).

\begin{lem} \label{Gk1012kTables} 
 \ The group $G_{k,1}(0,1|2|\ldots|k{\rm -}1)$ consists of the elements 
of $G_{k,1}$ with tables of the form
\[ \left[ \begin{array}{lllll lllll lll}
u_1 \ \dots \ u_{\ell} &|& p_1 a_2 \ \ \, \dots \ p_{\ell-1} a_2 &|& \dots
  \ \ \dots \ \ \dots &|& p_1 a_{k-1} \ \ \ \ \ \dots 
 \ p_{\ell-1} a_{k-1} \\   
v_1 \ \dots \ v_{\ell} &|& q_1^{(2)} a_2 \ \dots \ q_{\ell-1}^{(2)} a_2 &|& 
 \dots \ \ \dots \ \ \dots &|& q_1^{(k-1)} a_{k-1} \ \dots 
 \ q_{\ell-1}^{(k-1)} a_{k-1} 
\end{array}        \right] , \]
or equivalently,

\medskip

\hspace{0.3in}
$\{(u_r,\, v_r) : 1 \le r \le \ell\} \ \ \cup \ \ $
$\bigcup_{i=2}^{k-1} \,  \{(p_s a_i, \ q_s^{(i)} a_i) : \, $
                    $1 \le s \le \ell-1\}$.

\bigskip

\noindent
Here $\{u_r : 1 \le r \le \ell\}$ and $\{v_r : 1 \le r \le \ell\}$ are
maximal prefix codes over $\{a_0,a_1\}$ of equal cardinality $\ell \ge 1$,
and \ $\{p_s : 1 \le s \le \ell-1\} =$ 
${\sf spref}(\{u_r: 1 \le r \le \ell\})$.
The map $u_r \mapsto v_r$ (for $1 \le r \le \ell$) is an arbitrary bijection.
For every $\, i = 2, \ldots, k{\rm -}1:$
 \ $\{q_1^{(i)}, \, \ldots \, , q_{\ell-1}^{(i)}\}$ $=$
${\sf spref}(\{v_1, \ldots, v_{\ell}\})$, and the map
$p_j a_i \mapsto q_j^{(i)} a_i\,$  (for $1 \le j \le \ell-1$)
is an arbitrary bijection. 
\end{lem}
{\sc Proof.}  This follows in a straightforward way from the definition of
$G_{k,1}(0,1|2|\ldots|k{\rm -}1)$ and Lemma \ref{MaxprefCodeAk}. 
The fact that for every $a_i$, the number of elements $p_j a_i \,$ (and 
$q_j^{(i)} a_i$) is $\ell -1$ follows from the fact that the set of strict 
prefixes of $\{u_1, \dots, u_{\ell}\}$ is the set of interior vertices of 
the prefix tree of $\{u_1, \dots, u_{\ell}\}$ (over the alphabet 
$\{a_0,a_1\}$); see \cite[Lemma 4.7]{BiCoNP}.  
 \ \ \ $\Box$

\begin{defn} \label{DEFpFix}
A function $g$ {\em partially fixes} a set $S \subseteq A^*$
 \ iff \ $g(x) = x\,$ for every 
$\, x \in S \,\cap\, {\rm Dom}(g) \,\cap\, {\rm Im}(g)$.
This is also called partial pointwise stabilization.
For a subgroup $G \subseteq G_{k,1}$, the {\em partial fixator} (in $G$) of 
$S$ is

\medskip

\hspace{1.0in}  ${\rm pFix}_G(S) \ = \ \{g \in G : \ $
$(\forall x \in S \,\cap\, {\rm Dom}(g) \,\cap\, {\rm Im}(g))$
$[\,g(x) = x\,]\,\}$.
\end{defn}
If the set $S$ is a {\em right ideal} of $A^*$ then ${\rm pFix}_G(S)$ is 
a group \cite[Lemma 4.1]{BiCoNP}.

\begin{lem} \label{pFixTable} \hspace{-0.07in}.

\noindent {\rm (1)} \ The group ${\rm pFix}_{G_{2,1}}(a_0 \{a_0,a_1\}^*)$ 
consist of the elements of $G_{2,1}$ that have a table of the form
\[ \left[ \begin{array}{l lll}
a_0 & a_1 u_1 & \dots & a_1 u_{\ell} \\
a_0 & a_1 v_1 & \dots & a_1 v_{\ell} 
\end{array}        \right] , \]
where $\{u_1, \,\ldots, u_{\ell}\}$ and $\{v_1, \, \ldots, v_{\ell}\}$ are
maximal prefix codes over $\,\{a_0,a_1\}$.

\medskip

\noindent {\rm (2)} \ The subgroups 
$\,{\rm pFix}_{G_{2,1}}(a_1 \{a_0,a_1\}^*)\,$ and 
$\,{\rm pFix}_{G_{2,1}}(a_0 \{a_0,a_1\}^*)\,$ are isomorphic to $G_{2,1}$.
\end{lem}
{\bf Proof.} (1) The form of the tables follows immediately from the 
definition of ${\rm pFix}_{G_{2,1}}(a_0 \{a_0,a_1\}^*)$. 

\smallskip

\noindent (2) We define an isomorphism 
$\,\theta: G_{2,1} \to {\rm Fix}_{G_{2,1}}(a_1 \{a_0,a_1\}^*) \,$ by 
\[ \left[ \begin{array}{ccc}
x_1 & \ldots & x_n \\
y_1 & \ldots & y_n
\end{array} \right] \ \ \ \longmapsto 
 \ \ \ \left[ \begin{array}{l lll}
a_0 & a_1 x_1 & \dots & a_1 x_{\ell} \\
a_0 & a_1 y_1 & \dots & a_1 y_{\ell}
\end{array}        \right].  \]
\noindent This map is obviously a bijection, and it is easy to check that 
it is a homomorphism.
 \ \ \ $\Box$

\begin{defn} \label{DEFdictorder} {\bf (dictionary order).}
 \ For an alphabet $A_k = \{a_0,a_1, \,\ldots, a_{k-1}\}$, totally ordered 
as $\,a_0 < a_1 < \,\ldots\, < a_{k-1}$, the {\em dictionary order} on 
$A_k^{\,*}$ is defined as follows. For any $\,u, v \in A_k^{\,*}:$ 

\smallskip

\noindent $u \le_{\rm dict} v$ \ \ \ iff 

\smallskip

\noindent {\small \rm (1)} \ $u \le_{\rm pref} v$, \ \  or 

\noindent {\small \rm (2)} \ $u \not\le_{\rm pref} v$, and
there exist $p, s, t \in A_k^{\, *}$ and $\alpha, \beta \in A_k$ such 
that $u = p \alpha s$, $v = p \beta t$, and $\alpha < \beta$.
\end{defn}
In case {\small \rm (2)}, $p$ is the longest common prefix of $u$ and 
$v$, $\, \alpha$ is the next letter after $p$ in $u$, and $\beta$ is 
the next letter after $p$ in $v$.
Since case {\small \rm (2)} rules out case {\small \rm (1)}, $p$ is strictly 
shorter than $u$ and $v$, so the letters $\alpha$ and $\beta$ exist, 
and $\alpha \ne \beta$.

\medskip

\noindent From now on we assume that $A_k$ is an {\em ordered alphabet}, as 
in Def.\ \ref{DEFdictorder}.

\begin{defn} \label{DEFrank} {\bf (rank function for {\boldmath 
$\le_{\rm dict}$}).}
 \ Let $P \subset A_k^{\,*}$ be a finite set, and let 
$(p_1, \,\ldots\, , p_{\ell})$ be the list of all the elements of $P$ in 
increasing dictionary order. Then the {\em rank} of $p_j$ in $P$ is 
 \ ${\sf rank}_P(p_j) = j-1$.  Equivalently, 
 \ ${\sf rank}_P(p_j) \, = \, |\{q \in P:\,  q <_{\rm dict} p_j\}|$.
\end{defn}

\noindent The following concept is crucial for embedding $G_{2,1}$ into 
$G_{k,1}(0,1|2|\ldots|k{\rm -}1)$.

\begin{defn} \label{LEM2neighbor} {\bf ($*a_i$-successor).} 
 \ Consider any $\, a_i \in \{a_2, \ldots, a_{k-1}\}$.
Let $P \subset \{a_0,a_1\}^*$ be any finite maximal prefix code with 
$|P| \ge 2$, and let $(p_1, \, \ldots, p_{\ell})$ be the list of all the
elements of $P$ in increasing dictionary order on $\{a_0,a_1\}^*$, where 
$\ell = |P|$.  

For every $p_j \in P \smallsetminus \{p_1\}$, the {\em $*a_i$-successor} 
$\, (p_j)'_i \,$ of $p_j$ is the element of 
$\,{\sf spref}(P) \, a_i$, defined as follows, assuming 
$\, (p_{j+1})'_i, \ \ldots \ , (p_{\ell})'_i \,$ have already been chosen:

\medskip

 \ \ \ $(p_j)'_i \ = \ \min\{\, x a_i \in {\sf spref}(P) \, a_i \ :$ 
$ \ p_j <_{\rm dict} x a_i \ $ {\rm and}
$ \ x a_i \not\in \{(p_{j+1})'_i, \ \ldots \ , (p_{\ell})'_i\} \,\}$,

\medskip

\noindent where $\min$ uses the dictionary order in $\{a_0,a_1,a_i\}^*\,$ 
(i.e., over the three-letter alphabet $\{a_0,a_1, a_i\}$).   
\end{defn}
In other words, $(p_j)'_i\,$ is the nearest right-neighbor of $p_j$ in 
$\, {\sf spref}(P) \, a_i$ \ $(\subset \{a_0,a_1,a_i\}^*) \,$ that has not
yet been associated with another $p_m$ for $m > j$. In the definition of 
$p_{\ell}$, $\,\{(p_{j+1})'_i, \, \ldots \, , (p_{\ell})'_i\}$ $=$ 
$\varnothing$ \ (when $j = \ell$).

\medskip

\noindent {\bf Remarks.}
For the concept of $*a_i$-successor, the three-letter alphabet 
$\{a_0,a_1,a_i\}\,$ (for a chosen $a_i$, $2 \le i < k$) must not be 
confused with the $k$-letter alphabet $A_k$ (except, of course, when $k=3$).

Note that $(p_1)'_i$ is not defined.  Indeed, if $P$ has $\ell$ elements,
${\sf spref}(P)$ has only $\ell{\rm -}1$ elements; all $p_j \in P$ with 
$j>1$ have a $*a_i$-successor, so there is no element left in 
${\sf spref}(P) \, a_i$ to be the successor of $p_1$; this is further 
clarified by Lemma \ref{LEM2succFormula}, which gives a simple formula for
the $*a_i$-successor.

\begin{lem} \label{LEM2succFormula} {\bf ($*a_i$-successor formula).} 
 \ Let $\, a_i \in \{a_2, \ldots, a_{k-1}\}$, and let 
$P \subset \{a_0,a_1\}^*$ be a finite maximal prefix code with $|P| \ge 2$. 
Then every element of $P \smallsetminus a_0^{\,*}$ can be written 
(uniquely) in the form $\, u a_1 a_0^{\,m} \,$ (where 
$u \in \{a_0,a_1\}^*$ and $m \ge 0$), and its $*a_i$-successor is 

\medskip

\hspace{0.7in}  $\, (u a_1 a_0^{\,m})'_i \ = \ u a_i$.

\medskip

\noindent The elements of $a_0^{\,*}$ have no $*a_i$-successor.

\smallskip

Conversely, for every $\, u a_i \in {\sf spref}(P) \, a_i \,$ we have: 
 \ $u a_i \,$ is the $*a_i$-successor of $\, u a_1 a_0^{\,m}$, where $m$ is 
the unique number such that 
 \ $u a_1 a_0^{\,m} \in P \, \cap \, u a_1 a_0^{\,*}$. Different elements 
of $P \smallsetminus a_0{\,^*}$ have different $*a_i$-successors.  
\end{lem}
{\sc Proof.} According to the definition of $*a_i$-successor, the first 
element $p_1 \in P$ has no $*a_i$-successor. In every maximal prefix code, 
$p_1 \in a_0^{\,*}$, hence elements of $a_0^{\,*}$ have no $*a_i$-successor.

Every element of $\,\{a_0,a_1\}^* \smallsetminus a_0{\,^*}\,$ contains an 
occurrence of $a_1$, and hence is of the form $u a_1 a_0^{\,m}$, for some 
$u \in \{a_0,a_1\}^*$ and $m \ge 0$. And $m$ is unique since $P$ is a prefix
code.   

Obviously, $u a_1 a_0^{\,m}$ $<_{\rm dict} u a_i \,$ in the dictionary 
order determined by $\, a_0 < a_1 < a_i$.  
Let $(p_1, \, \ldots, p_{\ell})$ be the list of all the elements of $P$ in 
increasing dictionary order, where $\ell = |P|$.  For 
$u a_1 a_0^{\,m} \in $ $P \smallsetminus a_0^{\,*}$, let us denote 
$\, {\sf rank}_P(u a_1 a_0^{\,m})\,$ by $r-1$, so $u a_1 a_0^{\,m} = p_r$.
We want to show that if there exists $ \, v \in {\sf spref}(P)\,$ 
such that $\, u a_1 a_0^{\,m} <_{\rm dict} v a_i <_{\rm dict} u a_i$, then 
$\, v a_i \in \{(p_{r+1})'_i, \, \ldots \, , (p_{\ell})'_i\}$.
This will imply that $u a_i$ is the minimum element in $\, {\sf spref}(P)$ 
$\smallsetminus$  $\{(p_{r+1})'_i, \,\ldots\, , (p_{\ell})'_i\}$,
satisfying $p_r <_{\rm dict} u a_i$; hence, $u a_i = (p_r)'_i\,$  by
Def.\ \ref{LEM2neighbor}.

\smallskip

The relations $v \in \{a_0,a_1\}^*$ and $\,u a_1 a_0^{\,m}$   $<_{\rm dict}$ 
$v a_i <_{\rm dict} u a_i\,$ imply that $v = u a_1 x$ for some 
$x \in \{a_0,a_1\}^*$. 

If $m = 0$ then $p_r = u a_1$ is a prefix of $v = u a_i x$, which contradicts 
the fact that $v \in {\sf spref}(P)$. Hence, $(u a_1)'_i = u a_i$. 
In particular, the case $m=0$  applies to $p_{\ell}$, since in a maximal 
prefix code the last element belongs to $a_1^*$, i.e., 
$p_{\ell} = a_1^{\,n}$ for some $n > 0$; hence 
$(p_{\ell})'_i = a_1^{\,n-1} a_i$;  so the $*a_i$-successor formula holds. 

Let us now assume by induction on decreasing $r$ (ranging from $\ell$ down
to $2$), that for all $j > r$: 

\smallskip

 \ \ \ \ \ if there exists $v_j \in {\sf spref}(P)$ such that 
$\,p_j = u_j a_1 a_0^{\,m_j} <_{\rm dict} v_j a_i <_{\rm dict} u_j a_i$, 

 \ \ \ \ \ then 
        $\,v_j a_i \in \{(p_{j+1})'_i, \ \ldots \ , (p_{\ell})'_i\}$.

\smallskip

\noindent
The inductive assumption holds when $r = \ell$ (since $v_{\ell}$ does not 
exist, by the case $m=0$).
So we can assume that $m > 0$, since we already proved that 
$(u a_1)' = u a_i$. 
We want to prove the inductive hypothesis for $r$, i.e.: If there exists 
$v_r \in {\sf spref}(P)$ such that
$\, p_r = u a_1 a_0^{\,m} <_{\rm dict} v_r a_i <_{\rm dict} u a_i$, then
$v_r a_i \in \{(p_{r+1})'_i, \, \ldots \, , (p_{\ell})'_i\}$. 

\medskip

\noindent Case where $\, u a_1 a_0^{\,m} <_{\rm dict} v_r :$

Then $\, p_r = u a_1 a_0^{\,m} <_{\rm dict} v_r <_{\rm dict} v_r z $
$<_{\rm dict}$ $v_r a_i <_{\rm dict} u a_i\,$ for every 
$z \in \{a_0,a_1\}^+$.  Since $P$ is maximal, $v_r$ is the prefix of some 
$\, p_j = v_r z_j \in P$, hence  $\, p_r = u a_1 a_0^{\,m}$  $<_{\rm dict}$ 
$p_j <_{\rm dict} v_r a_i <_{\rm dict} u a_i$, hence $j>r$.  
Now by induction (since $j>r$), 
$\, v_r a_i \in \{(p_j)'_i, \, \ldots \, , (p_{\ell})'_i\} \,$ 
($\subseteq \{(p_{r+1})'_i, \, \ldots \, , (p_{\ell})'_i\}$).

\medskip

\noindent Case where $\, v_r = u a_1 x \le_{\rm dict} u a_1 a_0^{\,m}:$

Then $v_r = u a_1 a_0^{\,k}$ for some $k < m$ (we cannot have $k = m$, since
$v_r \in {\sf spref}(P)$, i.e., $v_r$ is a strict prefix).  In that case, 
$p_r = u a_1 a_0^{\,m} <_{\rm dict} v_r a_1 = u a_1 a_0^{\,k} a_1$ 
$\le_{\rm dict} u a_1 a_0^{\,k} a_1 z$
$<_{\rm dict} v_r a_i <_{\rm dict} u a_i$ for all $z \in \{a_0,a_1\}^*$. 
Since $P$ is maximal, $v_r a_1$ is a prefix of some 
$p_j = v_r a_1 z_j \in P$, hence 
$\,p_r <_{\rm dict} v_r a_1 \le_{\rm dict} p_j <_{\rm dict} v_r a_i$ 
$<_{\rm dict}$ $u a_i$, hence $j>r$.
Then by induction, $v_r a_i \in \{(p_j)'_i, \,\ldots\, , (p_{\ell})'_i\}\,$
($\subseteq \{(p_{r+1})'_i, \,\ldots\, , (p_{\ell})'_i\}$).

\medskip

In conclusion, if there exists $v_r \in {\sf spref}(P)$ such that
$\, p_r = u a_1 a_0^{\,m} <_{\rm dict} v_r a_i <_{\rm dict} u a_i$, then
$v_r a_i \in \{(p_{r+1})'_i, \,\ldots\, , (p_{\ell})'_i\}$. 
So, $(p_r)'_i = u a_i$. 

\medskip

For the converse, we saw in Lemma \ref{OmegaPref} that if $P$ is a 
finite maximal prefix code then every $w \in \{a_0,a_1\}^{\omega}$ has a 
unique prefix in $P$. By taking $w = u a_1 a_0^{\, \omega}$ we conclude that 
$P \cap u a_1 a_0^{\,*}$ is a singleton; i.e., $u$ determines $m$. Since $u$ 
determines $m$, it follows that the function 
$ \ u a_1 a_0^{\,m} \in P \smallsetminus a_0^{\,*} \,\longmapsto\, u a_i \ $ 
is injective.   The converse follows immediately now from the fact that 
$(u a_1 a_0^{\,m})'_i = u a_i$.
 \ \ \ $\Box$

\begin{lem} \label{LEM2succExists} 
 \ Let $a_i \in \{a_2, \ldots, a_{k-1}\}$, and let $P \subset \{a_0,a_1\}^*$ 
be a finite maximal prefix code, ordered as 
$\,p_1 <_{\rm dict} \, \dots \, <_{\rm dict} p_{\ell}$, where 
$\ell = |P| \ge 2$.  Then

\medskip

\noindent {\small \bf (1)} \ $P \, \cup \, {\sf spref}(P) \, a_i$ \ is a 
finite maximal prefix code over the three-letter alphabet $\{a_0,a_1,a_i\}$.

Moreover, $\, P \, \cup \, {\sf spref}(P) \, \{a_2, \ldots, a_{k-1}\}$ 
 \ is a finite maximal prefix code over $A_k$. 

\smallskip

\noindent {\small \bf (2)} \ $\{(p_2)'_i, \,\ldots\, , (p_{\ell})'_i\}$ 
$\,=\,$ ${\sf spref}(P) \,a_i$ \ ($\subset \{a_0,a_1,a_i\}^*$). 

\smallskip

\noindent {\small \bf (3)} \ Consider a one-step restriction, in which $P$ is
replaced by $\,P_r = (P \smallsetminus \{p_r\}) \, \cup \, p_r \{a_0,a_1\}$.
Then $\, (p_r a_0)'_i$ and $(p_r a_1)'_i$ are uniquely determined by $p_r$ 
as follows: 

\smallskip

$\bullet$
 \ if $2 \le r \le \ell$ then $\,(p_r a_1)'_i =  p_r a_i\,$ and 
$\, (p_r a_0)'_i = (p_r)'_i \,$;

\smallskip

$\bullet$ 
 \ if $r = 1$ then $p_1 \in P \cap a_0^{\,*}$, hence 
$p_1 a_0 \in P_1 \cap a_0^{\,*}$, and 
$p_1 a_1 \in P_1 \smallsetminus a_0^{\,*}$; so $(p_1)'_i$ and 
$(p_1 a_0)'_i$ do not exist; \ but $\,(p_1 a_1)'_i = p_1 a_i$ exists.
\end{lem}
{\sc Proof.} (1) is straightforward. (2) follows immediately from Lemma 
\ref{LEM2succFormula} and the fact that for any finite prefix code $P$ over
$\{a_0,a_1\}$, \ $|{\sf spref}(P)| = |P| - 1$.  

(3) Let $p_r = u a_1 a_0^{\,m} \in P \smallsetminus a_0^{\,m}$, so 
$(p_r)'_i = u a_i$ (by Lemma \ref{LEM2succFormula}).
By Lemma \ref{extStep}, $P_r$ is a maximal prefix code over $\{a_0,a_1\}$.
Applying Lemma \ref{LEM2succFormula} to $P_r$ and its elements 
$p_r a_0 = u a_1 a_0^{\,m+1}$ and $\,p_r a_1 = u a_1 a_0^{\,m} a_1$, we 
obtain $(p_r a_0)' = u a_i = (p_r)'_i$ and 
$(p_r a_1)'_i = u a_1 a_0^{\,m} a_i = p_r a_i$.

If $p_1 = a_0^{\,m} \in P \cap a_0^{\,*}$ (for some $m \ge 0$), then 
$p_1 a_0 = a_0^{\,m+1}$, and $p_1 a_1 = a_0^{\,m} a_1$. Hence $(p_1)'_i$ and
$(p_1 a_0)'_i$ do not exist.  But $(p_1 a_1)'_i = a_0^{\,m} a_i = p_1 a_i$.
 \ \ \ $\Box$

\begin{lem} \label{EmbedG21inGk1}
 \ For every $k \ge 3$ there exists an embedding 
$\, \iota: G_{2,1} \hookrightarrow G_{k,1}(0,1|2|\ldots|k{\rm -}1)$.
\end{lem}
{\sc Proof.} We define the embedding 
\[ g =  \left[ \begin{array}{lll}
p_1 & \ldots & p_{\ell}  \\
q_1 & \ldots & q_{\ell}
\end{array} \right] \ \ \ \stackrel{\iota}{\hookrightarrow}
\hspace{4.9in}  \]
\[ \hspace{0.in} \iota(g) = \left[ \begin{array}{ll lll lll ll ll llll}
a_0\! &|& a_1p_1 \ \ldots \ a_1p_{\ell} &| 
  &(a_1p_1)'_2 \ \ldots \ (a_1p_{\ell})'_2 
  &| & \ldots \ & \ \ldots &| 
  &(a_1p_1)'_{k-1} \ \ldots \ (a_1p_{\ell})'_{k-1}  \\
a_0\! &| &a_1q_1 \ \ldots \ a_1q_{\ell} &| 
  &(a_1q_1)'_2 \ \ldots \ (a_1q_{\ell})'_2 
  &| & \ldots \ & \ \ldots &| 
  &(a_1q_1)'_{k-1} \ \ldots \ (a_1q_{\ell})'_{k-1}
\end{array} \right], \]

\noindent
where $\{p_1, \ldots, p_{\ell}\}$ and $\{q_1, \ldots, q_{\ell}\}$ are finite
maximal prefix codes over $\{a_0, a_1\}$.
Equivalently, the table for $\iota(g)$ is

\smallskip

 \ \ \   \ \ \ $\{(a_0, \, a_0)\}$   $ \ \ \cup \ \ $ 
$\{(a_1p_r, \, a_1q_r) : 1 \le r \le \ell\}$     $ \ \ \cup \ \ $ 
$\bigcup_{i=2}^{k-1}$
$\big\{\big((a_1p_r)'_i, \ (a_1q_r)'_i\big) : 1 \le r \le \ell \big\}$.  

\bigskip

\noindent 
The function $\iota$ is well defined, as a function between tables. 
Indeed, for all $r = 1, \ldots, \ell:$ \ $(a_1 p_r)'$ determines $p_r$ (by 
Lemma \ref{LEM2succFormula}), which in turn determines $q_r$ (via the table
for $g$), which determines $(a_1 q_r)'$.  And $\iota$ is obviously 
injective. 

\medskip

To show that $\iota$ is also a map from $G_{2,1}$ to
$G_{k,1}(0,1|2|\ldots|k{\rm -}1)$, we show that the operation of 
{\em one-step restriction commutes with} $\iota$. 
Moreover, after that we can restrict tables so that when we compose two 
tables, the image row of the first table is equal to the domain row of the 
second; this makes it easy to show that $\iota$ is a homomorphism.

\medskip

For any $g \in G_{2,1}$, given by a table 
$\{(p_r, \,q_r) : 1 \le r \le \ell\}$, the restriction of $g$ at $p_s$ 
(for $1 \le s \le \ell$) is given by the table  

\smallskip

 \ \ \  \ \ \  ${\sf restr}_{p_s}(g)$  $ \ = \ $
$\{(p_s a_0, \ q_s a_0), \ (p_s a_1, \ q_s a_1)\}$
$ \ \cup \ $  $\{(p_j,\, q_j) : 1 \le j \le \ell, \ j \ne s\}$. 

\smallskip

\noindent Similarly, for any $f \in G_{k,1}$, given by a table
$\{(u_j, v_j) : 1 \le j \le m\}$, the restriction at $u_t$ (for 
$1 \le t \le m$) is given by the table

\smallskip

 \ \ \  \ \ \  ${\sf restr}_{u_t}(f)$ $ \ = \ $
$\{(u_t a_i, \, v_t a_i) : a_i \in A_k\}$
$ \ \cup \ $  $\{(u_j, v_j) : 1 \le j \le m, \ j \ne t\}$. 

\bigskip

\noindent (1) (Commutation:) \ Verification that 
 \ $\iota({\sf restr}_{p_r}(g)) = {\sf restr}_{a_1p_r}(\iota(g))$
 \ \ (for $r = 1, ..., \ell$):
\[ {\sf restr}_{p_r}(g) \ = \ \left[ \begin{array}{llllll}
\dots & p_{r-1} & p_r a_0 & p_r a_1 & p_{r+1} & \dots  \\
\dots & q_{r-1} & q_r a_0 & q_r a_1 & q_{r+1} & \dots
\end{array} \right] \ \ \ \stackrel{\iota}{\hookrightarrow} 
 \ \ \ \iota({\sf restr}_{p_r}(g)) \hspace{2.3in} \]
\[ = \left[ \begin{array}{lllll lllll ll}
\dots & a_1p_{r-1} & a_1p_ra_0 & a_1p_ra_1 & a_1p_{r+1} & \dots & \dots &
    (a_1p_{r-1})'_i & (a_1p_ra_0)'_i & (a_1p_ra_1)'_i & (a_1p_{r+1})'_i & \dots \\   
\dots & a_1q_{r-1} & a_1q_ra_0 & a_1q_ra_1 & a_1q_{r+1} & \dots & \dots &
    (a_1q_{r-1})'_i & (a_1q_ra_0)'_i & (a_1q_ra_1)'_i & (a_1q_{r+1})'_i & \dots
\end{array} \right] \hspace{0.in} \]
\[ = \left[ \begin{array}{lllll lllll ll}
\dots & a_1p_{r-1} & a_1p_ra_0 & a_1p_ra_1 & a_1p_{r+1} & \dots & \dots &
    (a_1p_{r-1})'_i & (a_1p_r)'_i & a_1p_r a_i & (a_1p_{r+1})'_i & \dots \\
\dots & a_1q_{r-1} & a_1q_ra_0 & a_1q_ra_1 & a_1q_{r+1} & \dots & \dots &
    (a_1q_{r-1})'_i & (a_1q_r)'_i & a_1q_r a_i & (a_1q_{r+1})'_i & \dots
\end{array} \right];  \hspace{0.1in} \]
the latter equality holds by Lemma \ref{LEM2succExists}. Here, $a_i$ 
ranges over $\{a_2, \ldots, a_{k-1}\}$. Recall that here the restriction 
${\sf restr}_{p_r}$ happens over the alphabet $\{a_0, a_1\}$.

\noindent 
Equivalently, the latter table for $\iota({\sf restr}_{p_r}(g))$ is 

\medskip

$\{(a_0,a_0)\} \ \ \cup \ \ $
$\{(a_1 p_r a_0, \ a_1 q_r a_0), \ (a_1 p_r a_1, \ a_1 q_r a_1)\}$
$  \ \ \cup \ \ $
$\{(a_1 p_j, \ a_1 q_j) : 1 \le j \le \ell, \ j \ne r\}$

\smallskip

$ \ \cup \ \ $     $\bigcup_{i=2}^{k-1}$
$\big( \{((a_1 p_r)'_i, \ (a_1 q_r)'_i), \ (a_1 p_r a_i, \ a_1 q_r a_i)\}$
$ \,\cup \,$
$\{((a_1 p_j)'_i, \ (a_1 q_j)'_i) : 1 \le j \le \ell, \ j \ne r\} \big)$.

\bigskip

\noindent On the other hand, the table for ${\sf restr}_{a_1p_r}(\iota(g))$
in $G_{k,1}$ is
\[ \left[ \begin{array}{lllll lllll llll lll}
\dots \ a_1p_{r-1} & a_1p_ra_0 & a_1p_ra_1 \ \ldots \ a_1p_r a_i \ \dots \  
    a_1p_{r+1} & \dots \ \dots &
    (a_1p_{r-1})'_i & (a_1p_r)'_i & (a_1p_{r+1})'_i & \dots \\
\dots \ a_1q_{r-1} & a_1q_ra_0 & a_1q_ra_1 \ \ldots \ a_1q_r a_i \ \dots \   
    a_1q_{r+1} & \dots \ \dots &
    (a_1q_{r-1})'_i & (a_1q_r)'_i & (a_1q_{r+1})'_i & \dots 
\end{array} \right]. \]
Recall that the restriction ${\sf restr}_{a_1 p_r}$ happens over the 
alphabet $A_k$. 
 
\noindent Equivalently, the latter table for 
${\sf restr}_{a_1p_r}(\iota(g))$ is

\medskip

$\{(a_0,a_0)\} \ $

\smallskip

$ \ \cup \ \ $
$\{(a_1 p_r a_0,\, a_1 q_r a_0), \ \ldots \ , (a_1 p_r a_i,\, a_1 q_r a_i),$
$ \ \ldots \ , (a_1 p_r a_{k-1},\, a_1 q_r a_{k-1}) \}$

\smallskip

$ \ \cup \ \ $
$\{(a_1 p_j, \ a_1 q_j) : 1 \le j \le \ell, \ j \ne r\}$

\smallskip

$ \ \cup \ \ $     $\bigcup_{i=2}^{k-1}$
$\{((a_1 p_j)'_i, \ (a_1 q_j)'_i) : 1 \le j \le \ell \}$.

\medskip

\noindent We see that the tables for $\iota({\sf restr}_{p_r}(g))$ and
${\sf restr}_{a_1 p_r}(\iota(g))$ are the same, up to the order of the
entries.  More precisely, 

\medskip

$\iota({\sf restr}_{p_r}(g))$ 

\smallskip

$= \ \{(a_0,a_0)\} \ \ \cup \ \ $
$\{(a_1 p_r a_0, \ a_1 q_r a_0), \ (a_1 p_r a_1, \ a_1 q_r a_1)\}$
$  \ \ \cup \ \ $
$\{(a_1 p_j, \ a_1 q_j) : 1 \le j \le \ell, \ j \ne r\}$

\smallskip

 \ \ \ \ \ $ \cup \ \ $     $\bigcup_{i=2}^{k-1}$
$\big( \{((a_1 p_r)'_i, \ (a_1 q_r)'_i), \ (a_1 p_r a_i, \ a_1 q_r a_i)\}$
$ \,\cup \,$
$\{((a_1 p_j)'_i, \ (a_1 q_j)'_i) : 1 \le j \le \ell, \ j \ne r\} \big)$.

\medskip

$= \ \{(a_0,a_0)\} \ \ \cup \ \ $
$\{(a_1 p_r a_0, \ a_1 q_r a_0), \ (a_1 p_r a_1, \ a_1 q_r a_1)\}$
$  \ \ \cup \ \ $
$\{(a_1 p_j, \ a_1 q_j) : 1 \le j \le \ell, \ j \ne r\}$

\smallskip

 \ \ \ \ \  $ \cup \ \ $     
$\{(a_1 p_r a_i, \ a_1 q_r a_i) :  2 \le i \le k-1\}$

\smallskip

 \ \ \ \ \ $ \cup \ \ $   $\bigcup_{i=2}^{k-1}$
$\big( \{((a_1 p_r)'_i, \ (a_1 q_r)'_i)\}$ 
$ \,\cup \,$
$\{((a_1 p_j)'_i, \ (a_1 q_j)'_i) : 1 \le j \le \ell, \ j \ne r\} \big)$.

\medskip

$= \ \{(a_0,a_0)\} \ $
$ \ \cup \ \ $
$\{(a_1 p_r a_0,\, a_1 q_r a_0), \ \ldots \ , (a_1 p_r a_i,\, a_1 q_r a_i),$
$ \ \ldots \ , (a_1 p_r a_{k-1},\, a_1 q_r a_{k-1}) \}$

\smallskip

 \ \ \ \ \ $ \cup \ \ $
$\{(a_1 p_j, \ a_1 q_j) : 1 \le j \le \ell, \ j \ne r\}$

\smallskip

 \ \ \ \ \ $ \cup \ \ $     $\bigcup_{i=2}^{k-1}$
$\{((a_1 p_j)'_i, \ (a_1 q_j)'_i) : 1 \le j \le \ell \}$

\smallskip

$= \ {\sf restr}_{a_1p_r}(\iota(g))$.

\bigskip

\noindent (2) To complete the proof that $\iota$ is a homomorphism, 
consider $h, g \in$ $G_{2,1}$ with tables
\[ g \ = \ \left[ \begin{array}{lll}
\dots & p_r & \dots \\
\dots & q_r & \dots
\end{array}        \right]
 \ \ \longmapsto \ \
\iota(g) \ = \ \left[ \begin{array}{lll llllll}
a_0 & \dots & a_1p_r & \dots & \dots & (a_1p_r)'_i & \dots & \dots \\
a_0 & \dots & a_1q_r & \dots & \dots & (a_1q_r)'_i & \dots & \dots
\end{array}        \right], \]
\[ h \ = \ \left[ \begin{array}{lll}
\dots & q_r & \dots \\
\dots & s_r & \dots
\end{array}        \right]
 \ \ \longmapsto \ \
\iota(h) \ = \ \left[ \begin{array}{lll llll ll}
a_0 & \dots & a_1q_r & \dots & \dots & (a_1q_r)'_i & \dots & \dots \\
a_0 & \dots & a_1s_r & \dots & \dots & (a_1s_r)'_i & \dots & \dots
\end{array}        \right] \, , \]
where $r = 1, \, \ldots \, , \ell$.

We can indeed assume that the output row of $\iota(g)$ is equal to the
input row of $\iota(h)$, since we proved that one-step restrictions commute 
with $\iota$. We just need the obtain row equality on the 
$\{a_0,a_1\}^*$-part of the table; the remainder of the rows are then equal
too, since $q_r$ determines $(a_1q_r)'_i$.  

Then by composing the tables we obtain 
\[ \iota(h) \circ \iota(g) \ = \ \left[ \begin{array}{lll llll ll}
a_0 & \dots & a_1p_r & \dots & \dots & (a_1p_r)'_i & \dots & \dots \\
a_0 & \dots & a_1s_r & \dots & \dots & (a_1s_r)'_i & \dots & \dots
\end{array}        \right] \, , \]
while $h \circ g$ has a table
\[ h \circ g \ = \ \left[ \begin{array}{lll}
\dots & p_r & \dots \\
\dots & s_r & \dots 
\end{array}        \right]
 \ \ \longmapsto \ \  
\iota(h \circ g) \ = \ \ \left[ \begin{array}{lll llll ll}
a_0 & \dots & a_1p_r & \dots & \dots & (a_1p_r)'_i & \dots & \dots \\
a_0 & \dots & a_1s_r & \dots & \dots & (a_1s_r)'_i & \dots & \dots
\end{array}        \right] \, . \]
So, $\iota(h) \circ \iota(g) = \iota(h \circ g)$.
 \ \ \   \ \ \ $\Box$

\bigskip

\noindent This completes the proof of Theorem \ref{THMembed2intok}.

\section{Comments}

{\bf (1)} Explanation of Higman's numbers $\, K = 1 + (k-1) \, d \,$ for 
$d \ge 1$:  

For $K > k$, Higman's embedding $G_{K,1} \le G_{k,1}$ uses a bijective
encoding of the alphabet $A_K = \{a_0, a_1, \, \ldots, a_{K-1}\}\,$ onto any 
maximal prefix code of size $K$ over $A_k$.  
Any maximal prefix code over the alphabet $A_k$ has cardinality 
$\,1 + (k-1) \, d \,$ for some $d \ge 1$, where $d$ is the number of 
interior vertices in the prefix tree of the maximal prefix code. 
This gives the possible values of $K$ for this method.

\medskip

\noindent {\bf (2)} For the embedding $G_{2,1} \le G_{k,1}$ with $k > 2$, 
the above encoding method does not work; obviously, there is no maximal 
prefix code of size 2 over $A_k$.

In \cite{BCS01}, $F_{2,1}$ is embedded into $F_{k,1}$ by ignoring the middle 
edges in every caret (just keeping the edges labeled by $a_0$ and $a_{k-1}$), 
where $a_{k-1}$ now represents the letter $a_1$ used by $F_{2,1}$. This works 
because (by the dictionary order preserving property of $F_{k,1}$), the 
mapping of the end edges of a caret determines the mapping of the
intermediary edges (see \cite{BCS01}). 
For $G_{k,1}$ with $k > 2$, this does not work: the mapping of the 
intermediary edges of a caret is not determined by the mapping of the end 
edges; for $G_{2,1}$ versus $G_{k,1}$, this ``embedding'' is not a function. 

\smallskip

Higman's embeddings, and the embedding $F_{h,1} < F_{k,1}$ for all 
$h, k \ge 2$, seem to be the basis of the general belief that all $G_{k,1}$ 
easily embed in one another. 
In Matte Bon's paper \cite{MatteBon}, the general embedding result is 
actually not stated explicitly, although it follows immediately from 
\cite[Cor.\ 11.16]{MatteBon}; however, in \cite{MatteBon} it is stated
more than once that these embeddings are ``well known''. 

\medskip

\noindent {\bf (3)} The idea for constructing an embedding 
$G_{2,1} \le G_{k,1}$ with $k > 2$ can be developed through a sequence of 
ideas that do not quite work. We illustrate this with $k = 3$, since this 
captures most of the difficulty. (This came about independently of Matte 
Bon's proof, which I did not know of.)

Consider any $G_{2,1}$-table $\,t = \{(p_j, q_j) : 1 \le j \le \ell\}$, 
where $P = \{p_j : 1 \le j \le \ell\}\,$ and 
$Q = \{q_j : 1 \le j \le \ell\}\,$
are finite maximal prefix codes over $A_2 = \{0,1\}$. 

\smallskip

\noindent $\bullet$ \ First idea: Map the $G_{2,1}$-table $t$ to the table
$\, \{(p_j, q_j) : 1 \le j \le \ell\} \,\cup\, \{(2,2)\}$.
This obviously  does not work, because $\,P \cup \{2\}\,$ is not a maximal 
prefix  code over $A_3 = \{0,1,2\}\,$ (except if $P = \{0,1\}$).

\smallskip

\noindent $\bullet$ \ Second idea: Use the fact that 
$\,P \,\cup\, {\sf spref}(P) \, 2\,$ is a maximal prefix  code over $A_3$ 
(Lemma \ref{MaxprefCodeAk}). 

Problem: How should $\,{\sf spref}(P) \, 2\,$ be mapped bijectively onto 
${\sf spref}(Q) \, 2 \ $?  \ The simple matching by dictionary order 
yields an injective map of $G_{2,1}$-tables to $G_{3,1}$-tables; but this 
map between tables is not a map from $G_{2,1}$ into $G_{3,1}$ 
(since the table map does not commute with restriction -- see the proof of
Lemma \ref{EmbedG21inGk1}).

\smallskip

\noindent $\bullet$ \ Third idea: Introduce the concept of $*2$-successor, 
and map $ \, t = \{(p_j, q_j) : 1 \le j \le \ell\}$ \ to
 \ $\{(p_j, q_j) : 1 \le j \le \ell\}$   $\,\cup\, $
$\{( (p_j)', (q_j)') : 2 \le j \le \ell\}$. 

Problem: Elements of $0^*$ have no $*2$-successor; we let $p_1 \in 0^*$,
which solves the problem in $P$ (we assume here that $P$ is written in
increasing dictionary order; but $Q$ isn't). But the $0^*$ element in $Q$ 
might not be $q_1$. So we still don't know how to map all the 
$*2$-successors.

\smallskip

\noindent $\bullet$ \ Fourth idea: Assume $P$ and $Q$ are ordered by the 
dictionary order of $A_2^{\,*}$ as 
$\,p_1 <_{\rm dict} \ \ldots \ <_{\rm dict} p_{\ell}$, and 
$\,q_1 <_{\rm dict} \ \ldots \ <_{\rm dict} q_{\ell}$.  
We map the table $t$ to the $G_{3,1}(0,1;2)$-table
\[ \ \ \  \ \ \  \left[ \begin{array}{lll lll}
p_1 &\dots & p_{\ell} & (p_2)' & \dots & (p_{\ell})' \\
q_{\pi(1)} & \dots & q_{\pi(\ell)} & (q_{\bar{\pi}(2)})' & \dots &
(q_{\bar{\pi}(\ell)})'
\end{array} \right], \hspace{1.in}  \]
where $\pi$ is a permutation of $\{1, 2, \ldots, \ell\}$, and $\bar{\pi}$ 
is the permutation of $\{2, \ldots, \ell\}$ determined by $\pi$ as follows.
For $i \in \{2, \ldots, \ell\}$,

\[ \bar{\pi}(i) \ = \
  \begin{cases}
  \pi(i) \ & \ \text{if \ $\pi(i) \in \{2, \,\ldots\, , \ell\}$ ,} \\
  \pi(1) \ & \ \text{otherwise (i.e., if $\, \pi(i) = 1$) .}
  \end{cases}       \hspace{1.5in}
\]
This is an injection of $G_{2,1}$ into $G_{k,1}$; but it is not a 
homomorphism because we do not have 
 \ $\ov{\rho \circ \pi} = \ov{\rho} \circ \ov{\pi}$ \ in general.
A counter-example is
\[ \pi \ = \ \left[ \begin{array}{llll}
1 & 2 & 3 & 4 \\
2 & 4 & 1 & 3
\end{array} \right], \ \
 \ \ \rho \ = \ \left[ \begin{array}{llll}
1 & 2 & 3 & 4 \\
3 & 4 & 1 & 2
\end{array} \right]. \]

\smallskip

\noindent $\bullet$ \ Finally: To the idea of the $*2$-successor, add the 
replacement of $G_{2,1}$ by its isomorphic copy 
$\,{\rm pFix}_{G_{2,1}}(0 \,\{0,1\}^*)$; this is carried out in section 2. 

\medskip

\noindent {\bf (4)} The asymmetry between the embeddings 
$G_{k,1} \le G_{2,1}$ and $G_{2,1} \le G_{k,1}$ is interesting. 
Matte Bon's \cite[Cor.\ 11.18]{MatteBon} implies that every copy of 
$G_{2,1}$ in $G_{k,1}$ (both acting on $A_k^{\, \omega}$) has a non-empty 
clopen set of global fixed points (in our Lemma \ref{EmbedG21inGk1}, this is
$a_0 A_k^{\, \omega}$), while such fixed points do not generally exist for 
$G_{k,1} \le G_{2,1}$.

Intuitively, the difficulty about embedding $G_{2,1}$ into $G_{k,1}$ comes 
from the fact that a $G_{k,1}$-table contains more entries than a 
$G_{2,1}$-table; how can a $G_{2,1}$-table determine this extra entries?
For $F_{2,1} \le F_{k,1}$, order-preservation makes that determination 
\cite{BCS01}. For $G_{k,1} \le G_{2,1}$ the problem is reversed, and coding 
solves the problem \cite{Hig74}.
For $G_{2,1} \le G_{k,1}$, the concept of $*a_i$-successor is the main idea 
for making a $G_{2,1}$-table determine a $G_{k,1}$-table; in addition, in
order to get a homomorphism, the idea of replacing $G_{2,1}$ by its 
isomorphic copy $\,{\rm pFix}_{G_{2,1}}(a_0 \{a_0,a_1\}^*)\,$ plays a crucial 
role.  This makes sense in view of \cite[Cor.\ 11.18]{MatteBon} (which I did 
not know when version 1 of this was written).

\bigskip

\bigskip

\noindent {\bf Acknowledgements:} 
I would like to thank Matt Brin for many discussions, in particular for 
pointing out the idea for the mutual embeddings of the groups $F_{k,1}$.

I'm grateful to Nicol\'as Matte Bon for contacting me after I posted the 
first version of this paper (arXiv.org/abs/1902.09414v1).

\bigskip



\bigskip 

{\small 
\noindent
Dept.\ of Computer Science, and CCIB \\  
Rutgers University -- Camden \\   
birget@camden.rutgers.edu
}

\end{document}